# Stabilité du comportement des marches aléatoires sur un groupe localement compact

## Driss Gretete

*LATP, Centre de Mathématiques et informatiques, université de Provence, Aix-Marseille I, 9, rue Frédéric Joliot-Curie, 13453 Marseille cedex 13, France. E-mail : gretete@cmi.univ-mrs.fr*

**Résumé.** Dans cet article nous démontrons un théorème de stabilité des probabilités de retour sur un groupe localement compact unimodulaire, séparable et compactement engendré. Nous démontrons que le comportement asymptotique de $F^{*(2n)}(e)$ ne dépend pas de la densité $F$ sous des hypothèses naturelles. A titre d'exemple nous établissons que la probabilité de retour sur une large classe de groupes résolubles se comporte comme $\exp(-n^{1/3})$.

**Abstract.** The main result of this paper is a theorem about the stability of the return probability of symmetric random walks on a locally compact group which is unimodular, separable and compactly generated. We show that the asymptotic behavior of $F^{*(2n)}(e)$ does not depend on the choice of the density $F$ under natural assumptions. As an example we show that the return probabilities for a large class of solvable groups of exponential volume growth, behaves like $\exp(-n^{1/3})$.

## 1. Introduction

Soit $G$ un groupe localement compact, compactement engendré et séparable. Nous démontrons que les probabilités de retour dans un voisinage de l'indentité d'une marche aléatoire définie par une densité symétrique bornée admettant un moment d'ordre 2 et minorée sur un voisinage ouvert de l'identité par une constante strictement positive, ne dépendent essentiellement que de la structure algébrique du groupe. Ce résultat généralise un énoncé de [9] concernant les groupes de type fini. Les outils mis en oeuvre sont la dominance spectrale et les traces dans les algèbres de von Neumann. Nous illustrons le résultat principal en calculant le comportement asymptotique de la probabilité de retour dans certains groupes résolubles unimodulaires. Ces exemples sont à rapprocher de [12, 13, 14] répondant à une question de [11]. Ces résultats ont été annoncés dans [7]. Nous donnons ici les preuves complètes.

## 2. Notations

On considère un groupe $G$ localement compact séparable, compactement engendré, unimodulaire. Soit $e$ l'élément unité de $G$, soit $K$ un voisinage compact symétrique de $e$, qui engendre $G$. Soit $\mu$ une mesure de Haar sur $G$. On définit la longueur d'un élément $g \in G$ relativement à $K$ par

$$|g|_K = \min\{n \in \mathbb{N} \mid g \in K^n\}.$$

On remarque que pour tout $x, y \in G, |xy|_K \leq |x|_K + |y|_K$. Ceci permet de définir une distance sur $G$ en posant pour tout $x, y \in G$,

$$d_K(x, y) = |x^{-1}y|_K.$$

Pour $x \in G$, on notera

$$B_K(x, r) = \{y \in G \mid d_K(x, y) \leq r\}.$$

On définit pour deux fonctions $f$ et $g$ à valeurs réelles et définies sur une partie $I$ non majorée de $\mathbb{R}$, les relations

$$f \preceq g \iff \exists a, b > 0, \ \exists r > 0, \ \forall t \in I \cap ]r, +\infty[; bt \in I \cap ]r, +\infty[ \text{ et } f(t) \leq ag(bt),$$
$$f(t) \simeq g(t) \iff (f(t) \preceq g(t) \text{ et } g(t) \preceq f(t)).$$

Dans tout ce travail, on entend par comportement asymptotique d'une fonction ou une suite, la donnée d'un représentant de la relation $\simeq$.

Pour toute la suite, pour toute partie $U$ de $G$, on note $1_U$ la fonction indicatrice de $U$. Soit $\nu$ une mesure de probabilité symétrique sur $G$. On considère $\Omega = G^{\mathbb{N}}$, muni de la structure borélienne produit, $P = \delta_e \otimes \nu^{\otimes \mathbb{N}}$ la probabilité produit sur $\Omega$, où $\delta_e$ désigne la mesure de Dirac au point $e$.

$X_n : \Omega \to G$ la projection canonique, $X_0$ la variable certaine égale à $e$.

$Z_n(\omega) = \prod_{i=0}^n X_i(\omega); \omega \in \Omega$ définit comme dans [9] la marche sur $G$ associée à $\nu$. Si $A \subset G$ est un borélien, alors $P(Z_n \in A)$ est la probabilité de rejoindre $A$ au temps $n$, en partant de $e$ à l'instant 0.

Notre objectif est d'étudier le comportement asymptotique de $P(Z_{2t} \in B)$, et d'établir que ce comportement asymptotique ne dépend pas du choix de la densité lorsque $\nu$ est une mesure de probabilité à densité $F$ vérifiant les propriétés suivantes :

(i) $F$ est dans $L^1(G) \cap L^\infty(G)$.
(ii) Il existe un ouvert $U$ relativement compact symétrique qui engendre $G$, voisinage de $e$, tel que $F$ est minorée sur $U$ par une constante strictement positive.
(iii) $F$ est symétrique.

Pour deux fonctions $F$ et $H$ éléments de $L^2(G)$, on définit le produit de convolution $F * H$ par

$$(F * H)(x) = \int_G F(xy^{-1}) H(y) \, d\mu(y),$$

ainsi, on définit $F^{*2} = F * F$, et pour tout entier naturel $p$, $F^{*(p+1)} = F * F^{*p}$. Soit $F$ une fonction mesurable sur $G$ à valeurs réelles, intégrable et bornée, on définit les opérateurs de convolution à gauche et à droite par $F$ notés respectivement $L_F$ et $R_F$ sur $L^2(G)$ par

$$L_F h = F * h, \qquad R_F h = h * F.$$

Pour une mesure $\nu$ sur $G$ positive bornée, on définit les opérateurs $L_\nu$ et $R_\nu$ par

$$L_\nu h(x) = \int_G h(y^{-1}x) \, d\nu(y)$$

et

$$R_\nu h(x) = \int_G h(xy^{-1}) \, \mathrm{d}\nu(y).$$

Remarquons que si $\nu$ est une mesure bornée à densité $F$ relativement à la mesure de Haar $\mu$ alors $R_F = R_\nu$ et $L_F = L_\nu$.

Notation : On note $I$ l'application identique de $L^2(G)$, l'opérateur $\Delta_F = I - R_F$ est le laplacien associé à $F$, plus généralement, pour une mesure de probabilité $\nu$, on définit le laplacien associé à $\nu$ par $\Delta_\nu = I - R_\nu$.

2.1. *Définition de la probabilité de retour*

**Proposition 1.** *Soit $G$ un groupe localement compact unimodulaire moyennable, à génération compacte, soit $F$ une densité de probabilité sur $G$ symétrique et bornée, on suppose qu'il existe un voisinage ouvert $U$ de $e$, symétrique, relativement compact qui engendre $G$ tel que*

$$\inf\{F(x), x \in U\} > 0.$$

*Alors, pour tout voisinage $V$ de $e$ relativement compact, on a :*

$$F^{*(2t)}(e) \simeq P(Z_{2t} \in V).$$

**Démonstration.** Posons

$$p_t(x, y) = F^{*t}(x^{-1}y)$$

on a

$$p_t(x, y) = \iint_G p_1(x, z_1) \cdots p_1(z_{t-1}, y) \, \mathrm{d}\mu(z_1) \cdots \mathrm{d}\mu(z_{t-1}) = p_t(y, x)$$

donc $F^{*t}$ est symétrique.

$$F^{*(2t)}(g) = \int_G F^{*t}(x) F^{*t}(x^{-1}g) \, \mathrm{d}\mu(x) \leq \left( \int_G F^{*t}(x)^2 \, \mathrm{d}\mu(x) \right)^{1/2} \left( \int_G F^{*t}(x)^2 \, \mathrm{d}\mu(x^{-1}g) \right)^{1/2} = F^{*(2t)}(e).$$

Donc en intégrant on obtient

$$P(Z_{2t} \in V) = \int_V F^{*(2t)}(g) \, \mathrm{d}\mu(g) \leq \mu(B) F^{*(2t)}(e). \qquad \square$$

Reste à montrer l'autre domination, pour cela on utilisera le lemme suivant.

**Lemme 1.** *Soit $\nu$ une mesure de probabilité sur $G$, à densité $F$ relativement à la mesure de Haar, symétrique et bornée. Alors*

$$\|L_\nu\|_2^2 = \lim_{n \to +\infty} (F^{*(2n)}(e))^{1/2n}$$

*et ainsi*

$$\lim_{n \to +\infty} \left( \frac{F^{*(2n+2m)}(e)}{F^{*(2n)}(e)} \right) = \|L_\nu\|_2^{2m}.$$

**Démonstration (Comparer avec [17], Lemma 10.1, p. 110).** Posons $a_n = F^{*(2n)}(e)$ et $p_n(x,y) = F^{*n}(x^{-1}y)$. On a $\int_G F^{*n}(x)\,d\mu(x) = 1$, et donc par application de l'inégalité de Cauchy–Schwarz $a_n = \int_G F^{*n}(x)^2\,d\mu(x) > 0$, et on a

$$p_{2n+k}(e,e) = \int_G p_{2n}(e,g)p_k(g,e)\,d\mu(g) \leq p_{2n}(e,e)\int_G p_k(g,e)\,d\mu(g) = p_{2n}(e,e),$$

on en déduit en particulier que $0 < \frac{a_{n+1}}{a_n} \leq 1$. D'autre part

$$p_{2n}(e,e) = \int_G p_{n+1}(e,g)p_{n-1}(g,e)\,d\mu(g)$$
$$\leq \left(\int_G p_{n+1}(e,g)^2\,d\mu(g)\right)^{1/2}\left(\int_G p_{n-1}(g,e)^2\,d\mu(g)\right)^{1/2} = p_{2n-2}(e,e)^{1/2}p_{2n+2}(e,e)^{1/2}$$

ainsi $a_n^2 \leq a_{n-1}a_{n+1}$, donc $(\frac{a_{n+1}}{a_n})$ est croissante et alors admet une limite $l$, par suite $(a_n^{1/n})$ converge vers $l$.

Montrons maintenant que $l \leq \|L_\nu\|_2^2$.

On a,

$$L_\nu F^{*(n-1)}(x) = \int_G F^{*(n-1)}(g^{-1}x)F(g)\,d\mu(g) = F^{*n}(x)$$

d'autre part

$$\|F^{*n}\|_2^2 = \int_G F^{*n}(g)^2\,d\mu(g) = F^{*(2n)}(e)$$

donc

$$\left\|L_\nu \frac{F^{*(n-1)}}{\|F^{*(n-1)}\|_2}\right\|_2^2 = \frac{F^{*(2n)}(e)}{F^{*(2n-2)}(e)},$$

par conséquent

$$\|L_\nu\|_2^2 \geq \frac{F^{*(2n)}(e)}{F^{*(2n-2)}(e)},$$

ceci montre que $\|L_\nu\|_2^2 \geq l$.

Pour montrer que $\|L_\nu\|_2^2 \leq l$, on utilise : $\nu^{2n}(K) \leq F^{*(2n)}(e)\mu(K)$ où $K$ est un voisinage compact de $e$.

En utilisant ([2], théorème 1, p. 174)

$$\limsup \nu^{2n}(K)^{1/2n} = \|L_\nu\|_2^2 \leq \lim_n (F^{*(2n)}(e)\mu(K))^{1/2n} = l,$$

ceci montre que $\frac{a_{n+1}}{a_n}$ converge vers $\|L_\nu\|_2^2$ par suite $\frac{a_{n+m}}{a_n} = \prod_{k=0}^{n-1} \frac{a_{n+k+1}}{a_{n+k}}$ converge vers $\|L_\nu\|_2^{2m}$.     □

Maintenant nous allons achever la preuve de la proposition 1.

**Démonstration de la proposition 1.** Le groupe $G$ est moyennable, donc en appliquant le lemme 1, on obtient

$$\lim_n \frac{F^{*(2n+2m)}(e)}{F^{*(2n)}(e)} = 1.$$

On peut donc choisir $n$ assez grand pour que

$$\frac{F^{*(2n+2m)}(e)}{F^{*(2n)}(e)} \geq \frac{1}{2}.$$

Puisque $V$ est relativement compact, pour toute constante fixée $C_0$, on peut choisir $m$ assez grand pour qu'on ait pour tout $x \in V$, l'inégalité $d_U(e,x) \leq \frac{m^{1/2}}{4C_0}$, donc pour $m$ et $n$ assez grand on a successivement

$$r(n,m) = \frac{m^{1/2}}{2C_0} \frac{F^{*(2n+2m)}(e)}{F^{*(2n)}(e)} \geq \frac{m^{1/2}}{2C_0} \geq d_U(e,x).$$

Par conséquent, en appliquant ([8], théorème 4.2, p. 683), on a l'existence d'une constante $C_0$ telle que pour $x \in G$ tel que l'inégalité

$$F^{*(2n+2m)}(x) \geq \frac{1}{2} F^{*(2n+2m)}(e), \tag{1}$$

donc en intégrant l'inégalité (1) sur $V$ on obtient la domination souhaitée. □

## 3. Une propriété des traces

Nous établissons d'abord une propriété de domination inspirée de [3] avec une légère modification.

**Proposition 2.** *Soit $\mathcal{A}$ une algèbre de von Neumann, $\tau$ une trace normale fidèle semi finie sur $\mathcal{A}$. Soient $x, y$ deux éléments de $\mathcal{A}$, autoadjoints positifs, tels que $x \leq y$ et $\mathrm{sp}(x)$, $\mathrm{sp}(y)$ inclus dans $[0,1]$, $g$ une fonction continue croissante sur $[0,1]$, et $r \in ]0,1[$. On suppose $g$ positive sur $[r,1]$. Alors la fonction $h = g 1_{[r,1]}$ vérifie $\tau(h(x)) \leq \tau(h(y))$.*

**Démonstration.** On va détailler la démonstration en plusieurs étapes.

Etape 1 : $g$ étant croissante et continue sur $[0,1]$, il existe donc une suite croissante $(g_n)$ de fonctions en escaliers croissantes, qui converge uniformément vers $g$, en effet, on considère pour tout entier naturel $n$, la fonction en escalier $g_n$ définie sur $[0,1]$ par $g_n(x) = g(\frac{k}{2^n})$ si $x \in [\frac{k}{2^n}, \frac{k+1}{2^n}[$, où $0 \leq k \leq 2^n - 1$ et $g_n(1) = g(1)$.

La suite $(g_n)$ est évidement croissante car si $x \in [\frac{k}{2^n}, \frac{k+1}{2^n}[$ et $x \in [\frac{k'}{2^{n+1}}, \frac{k'+1}{2^{n+1}}[$ alors $\frac{k}{2^n} \leq \frac{k'}{2^{n+1}}$ par suite $g_n(x) \leq g_{n+1}(x)$, la croissance des $g_n$ découle immédiatement de celle de $g$.

Soit maintenant $\epsilon > 0$. La continuité uniforme de $g$ assure

$$\exists \delta > 0, \ \forall u, v \in [0,1]; \quad |u-v| \leq \delta \implies |g(u) - g(v)| \leq \epsilon.$$

Soit $N$ un entier tel que pour tout $n \geq N$, $\frac{1}{2^n} \leq \delta$. Pour $x \in [0,1[$, il existe un entier $k > 0$ tel que $x \in [\frac{k}{2^n}, \frac{k+1}{2^n}[$. On a

$$\left| x - \frac{k+1}{2^n} \right| \leq \delta,$$

donc

$$|g(x) - g_n(x)| \leq \epsilon.$$

Etape 2 : Pour tout entier $n > 0$, la fonction $h_n = 1_{[r,1]} g_n$ peut s'écrire comme combinaison linéaire à coefficients positifs de fonctions de la forme $1_{[s,1]}$, où $s > 0$.

En effet : on écrit $g_n = \sum_{k=0}^{m-1} \lambda_k 1_{[s_k, s_{k+1})}$, où $(s_k)_{0 \leq k \leq m}$ une subdivision de $[0,1]$ adaptée à $g_n$. Par définition le symbole $1_{[s_k, s_{k+1})}$ est égal à $1_{[s_k, s_{k+1}[}$ pour tout $k \in 0, \ldots, m-2$ et égal à $1_{[s_k, s_{k+1}]}$ pour $k = m-1$. Quitte à remplacer $(s_k)_{0 \leq k \leq m}$ par une subdivision plus fine, on pourra supposer qu'il existe $j \in \{1, \ldots, m-1\}$ tel que $s_j = r$.

Alors $h_n = \sum_{k=0}^{m-1} \lambda_k 1_{[s_k,s_{k+1})} 1_{[r,1]} = \sum_{k=j}^{m-1} \lambda_k 1_{[s_k,s_{k+1})} 1_{[r,1]} = \sum_{k=j}^{m-1} \lambda_k 1_{[s_k,s_{k+1})}$ alors

$$h_n = \sum_{k=j}^{m-1} \lambda_k 1_{[s_k,s_{k+1})} = \sum_{k=j}^{m-2} \lambda_k (1_{[s_k,1]} - 1_{[1,s_{k+1}]}) + \lambda_{m-1} 1_{[s_{m-1},1]} = \lambda_j 1_{[r,1]} + \sum_{k=j+1}^{m-1} (\lambda_k - \lambda_{k-1}) 1_{[s_k,1]}.$$

Etape 3 : D'après ([3], p. 170), pour tout $s > 0$, $\tau(1_{[s,+\infty[}(x)) \leq \tau(1_{[s,+\infty[}(y))$, et puisque le spectre de chacun des opérateurs est inclus dans $[0,1]$, il vient que

$$\tau(1_{[s,1]}(x)) \leq \tau(1_{[s,1]}(y)).$$

Etape 4 : On a $h(x) = \sup\{h_n(x), n \in \mathbb{N}\}$, et de même $h(y) = \sup\{h_n(y), n \in \mathbb{N}\}$, en effet pour tout $n \in \mathbb{N}$, $0 \leq h_n(x) \leq h(x) \leq \|g\|_\infty$, donc $h_n$ est bornée, d'où par application du théorème spectral, $h_n(x)$ converge fortement vers $h(x)$, par conséquent pour tout vecteur $\xi \in H$, $\lim_n (h_n(x)\xi|\xi) = (h(x)\xi|\xi)$.

Soit $u$ un opérateur de $L^2(G)$ tel que $\forall n, h_n(x) \leq u$ alors

$$\forall \xi \in L^2(G), \quad (h_n(x)\xi|\xi) \leq (u\xi|\xi).$$

Par passage à la limite on trouve $(h(x)\xi|\xi) \leq (u\xi|\xi)$, donc $h(x) \leq u$. Par suite $h(x) = \sup\{h_n(x), n \in \mathbb{N}\}$.

Etape finale : $h_n(x)$ et $h_n(y)$ sont des combinaisons linéaires à coefficients positifs de fonctions de la forme $1_{[s,1]}$ où $s > 0$, donc d'après le résultat de l'étape 3, on a l'inégalité $\tau(h_n(x)) \leq \tau(h_n(y))$, comme $h(x) = \sup_n h_n(x)$, par normalité de la trace on obtient

$$\tau(h(x)) = \sup_n \tau(h_n(x)) \leq \sup_n \tau(h_n(y)) = \tau(h(y)),$$

d'où l'inégalité. $\square$

**Théorème 1.** *Soit $\mathcal{A}$ une algèbre de von Neumann, $\tau$ une trace normale fidèle semi finie sur $\mathcal{A}^+$. Soient $x$ et $y$ deux éléments autoadjoints positifs de $\mathcal{A}$ de spectres compris dans $[0,1]$ tels qu'il existe une constante $C > 0$ vérifiant $I - x \leq C(I - y)$ alors*

$$\tau(y^n) \preceq \tau(x^n) + e^{-n}.$$

Pour démontrer ce théorème on a besoin du lemme suivant

**Lemme 2.** *Soit $c$ un réel dans $]0,1[$ et $r = 1 - c$. Pour $t \geq 1$, on a les inégalités*

  (i) $\forall \lambda \in [0,r]$, $\lambda^{2t} \leq \lambda e^{-ct}$.
  (ii) $\forall \lambda \in [r,1]$, $\lambda^{2t} \leq e^{(\lambda-1)t} - e^{-ct} + \lambda e^{-ct}$.
  (iii) $\forall \lambda \in [0,1]$, $e^{-2t}(e^{2\lambda t} - 1) \leq 2\lambda t e^{-t} + \lambda^t$.

**Démonstration.**

  (i) pour $\lambda \in ]0,r]$, $\lambda^{2t} = \lambda^t e^{t\ln\lambda} \leq \lambda e^{t(\lambda-1)} \leq \lambda e^{-ct}$.
  (ii) pour $\lambda \in [r,1]$, $\lambda^{2t} = \lambda^t e^{t\ln\lambda} \leq \lambda e^{t(\lambda-1)} = \lambda(e^{t(\lambda-1)} - e^{-ct}) + \lambda e^{-ct} \leq e^{t(\lambda-1)} - e^{-ct} + \lambda e^{-ct}$.
  (iii) pour tout $\lambda \in [0,1]$,

$$e^{\lambda 2t} - 1 = \sum_{k=1}^{+\infty} \frac{(\lambda 2t)^k}{k!} \leq \lambda 2t \sum_{k=0}^{+\infty} \frac{(\lambda 2t)^k}{(k+1)!} \leq 2\lambda t e^{2\lambda t}$$

donc

$$e^{-2t}(e^{2\lambda t} - 1) \leq 2\lambda t e^{(\lambda-1)2t}$$

par suite pour $\lambda \in [0, 1/2]$,

$$\mathrm{e}^{-2t}(\mathrm{e}^{2\lambda t} - 1) \leq 2\lambda t \mathrm{e}^{-t}$$

d'autre part pour $\lambda \in [1/2, 1]$,

$$2\lambda - 2 \leq \ln \lambda$$

et dans ce cas

$$\mathrm{e}^{-2t}(\mathrm{e}^{2\lambda t} - 1) = \mathrm{e}^{(\lambda-1)2t} - \mathrm{e}^{-2t} \leq \mathrm{e}^{(\lambda-1)2t} \leq \mathrm{e}^{t\ln\lambda} = \lambda^t$$

d'où l'inégalité. □

**Démonstration du théorème 1.** On prend $c = C^{-1}$ en choisissant $C$ assez grand pour que $c \in \,]0,1[$.

Par le théorème spectral, il existe des résolutions spectrales telles que $y = \int_{[0,1]} \lambda \, \mathrm{d}E_\lambda$ et $x = \int_{[0,1]} \lambda \, \mathrm{d}E'_\lambda$. Soit $t \geq 1$, on a d'après ([6], propositon 2, p. 3), si $T$ est dans l'algèbre de von Neumann $\mathcal{A}$, alors pour tout borélien $\Omega$, l'opérateur $1_\Omega(T)$ est aussi dans $\mathcal{A}$, on peut donc écrire

$$\tau(y^{2t}) = \tau\left(\int_{[0,1]} \lambda^{2t} \, \mathrm{d}E_\lambda\right) = \tau\left(\int_{[0,r[} \lambda^{2t} \, \mathrm{d}E_\lambda\right) + \tau\left(\int_{[r,1]} \lambda^{2t} \, \mathrm{d}E_\lambda\right),$$

et en utilisant le point (i) du lemme 2,

$$\tau(y^{2t}) \leq \tau\left(\int_{[0,r[} \lambda \mathrm{e}^{-ct} \, \mathrm{d}E_\lambda\right) + \tau\left(\int_{[r,1]} \lambda^{2t} \, \mathrm{d}E_\lambda\right) \leq \mathrm{e}^{-ct}\tau(y) + \tau\left(\int_{[r,1]} \lambda^{2t} \, \mathrm{d}E_\lambda\right)$$

en utilisant le point (ii) du lemme 2,

$$\tau\left(\int_{[r,1]} \lambda^{2t} \, \mathrm{d}E_\lambda\right) \leq \tau\left(\int_{[r,1]} (\mathrm{e}^{(\lambda-1)t} - \mathrm{e}^{-ct}) \, \mathrm{d}E_\lambda\right) + \tau\left(\int_{[r,1]} \lambda \mathrm{e}^{-ct} \, \mathrm{d}E_\lambda\right),$$

par conséquent

$$\tau(y^{2t}) \leq 2\mathrm{e}^{-ct}\tau(y) + \tau\left(\int_{[r,1]} (\mathrm{e}^{(\lambda-1)t} - \mathrm{e}^{-ct}) \, \mathrm{d}E_\lambda\right).$$

Soit $h$ la fonction définie sur $[0, 1]$ par $h(\lambda) = 1_{[r,1]}(\lambda)(\mathrm{e}^{(\lambda-1)t} - \mathrm{e}^{-ct})$, $h$ est positive car pour $\lambda \geq r$, on a $\lambda - 1 \geq r - 1 = -c$. L'inégalité ci-dessus s'écrit donc

$$\tau(y^{2t}) \leq 2\tau(y)\mathrm{e}^{-ct} + \tau(h(y)).$$

La trace de $I - cI + cx$ est infinie, on tronque cet opérateur à l'aide de la fonction $h$, puisque $y \leq I - cI + cx$, alors d'après la proposition 2, on obtient

$$\tau(h(y)) \leq \tau(h(I - cI + cx)).$$

D'autre part

$$h(I - cI + cx) = \int_{[0,1]} h(1 - c + c\lambda) \, \mathrm{d}E'_\lambda \leq \int_{[0,1]} (\mathrm{e}^{(-c+c\lambda)t} - \mathrm{e}^{-ct}) \, \mathrm{d}E'_\lambda$$

et en appliquant le (iii) du lemme ci-dessus, en remplaçant $t$ par $ct/2$ et en supposant $t > 2/c$ on obtient

$$\forall \lambda \in [0, 1], \quad \mathrm{e}^{-ct}(\mathrm{e}^{\lambda ct} - 1) \leq \lambda ct \mathrm{e}^{-ct/2} + \lambda^{ct/2}$$

donc

$$\tau(h(I-cI+cx)) \leq \tau\left(\int_{[0,1]} (\lambda c t e^{-ct/2} + \lambda^{ct/2})\, dE'_\lambda\right) = c t e^{-ct/2}\tau(x) + \tau\left(\int_{[0,1]} \lambda^{ct/2}\, dE'_\lambda\right)$$

d'où

$$\tau(y^{2t}) \preceq e^{-t} + \tau(x^{2t}). \qquad \square$$

## 4. Invariance des probabilités de retour

### 4.1. Comparaison des formes de Dirichlet

On considère toujours un groupe localement compact séparable unimodulaire, et $\nu$ une mesure de probabilité symétrique, on lui associe le laplacien défini par $\Delta_\nu = I - R_\nu$.

La forme de Dirichlet associée au laplacien est la forme quadratique définie sur $L^2(G)$ par

$$\varepsilon_\nu(f,f) = (\Delta_\nu f | f).$$

**Proposition 3.** *Si $\nu$ est une mesure de probabilité sur $G$ à densité $F$ symétrique bornée, alors :*

$$\varepsilon_\nu(f,f) = \frac{1}{2} \int_G \int_G |f(x) - f(xy)|^2 F(y)\, d\mu(x)\, d\mu(y).$$

**Démonstration (comparer avec [15], p. 97).** On a

$$\varepsilon_\nu(f,f) = \int_G \int_G \overline{f(x)}(f(x) - f(xy^{-1})) F(y)\, d\mu(x)\, d\mu(y)$$

par changement de variable, on obtient

$$\varepsilon_\nu(f,f) = \int_G \int_G \overline{f(x)}(f(x) - f(z)) F(z^{-1}x)\, d\mu(x)\, d\mu(z).$$

En échangent les rôles des variables $x$ et $z$ et tenant compte de la symétrie de $F$, on trouve

$$\varepsilon_\nu(f,f) = \int_G \int_G \overline{f(z)}(f(z) - f(x)) F(z^{-1}x)\, d\mu(x)\, d\mu(z)$$

en ajoutant ces dernières égalités on obtient

$$\varepsilon_\nu(f,f) = \frac{1}{2} \int_G \int_G |f(x) - f(z)|^2 F(z^{-1}x)\, d\mu(x)\, d\mu(z).$$

En faisant à nouveau le changement de variable $x^{-1}z = y$, on obtient le résultat. $\qquad \square$

Maintenant nous allons établir un résultat analogue à celui de [9], il s'agit de la comparaison des formes de Dirichlet qui va nous permettre par des raisonnement sur les traces de comparer les comportements des marches aléatoires associées à des mesures différentes vérifiant certaines hypothèses.

**Proposition 4.** *Soit $G$ un groupe localement compact séparable, compactement engendré unimodulaire, $\nu_1$ et $\nu_2$ des mesures de probabilité à densités $F_1, F_2$, relativement à la mesure de Haar sur $G$, telles que $F_1$ et $F_2$ soient symétriques intégrables et bornées sur $G$. On suppose qu'il existe un ouvert $U$ relativement compact voisinage de $e$ qui engendre $G$ tel que $\inf_{U^2} F_2 \geq r > 0$, et que $F_1$ admet un moment d'ordre 2 relativement à $U$, c'est à dire que $\int_G |x|_U^2 F_1(x)\, d\mu(x)$ est finie. Alors il existe une constante $C > 1$ telle que $\varepsilon_{\nu_1} \leq C\varepsilon_{\nu_2}$.*

**Démonstration.** Dans cette preuve on s'inspire de ([8], p. 682) et de ([15], Proposition VII.3.2).

On a pour $f \in L^2(G)$ pour $x \in G, y \in U$, pour tout $z \in U$,

$$|f(x) - f(xy)| \leq |f(x) - f(xz)| + |f(xz) - f(xy)|,$$

donc

$$|f(x) - f(xy)|^2 \leq 2|f(x) - f(xz)|^2 + 2|f(xz) - f(xy)|^2,$$

et en intégrant on obtient

$$|f(x) - f(xy)|^2 \leq 2\mu(U)^{-1} \int_U |f(x) - f(xz)|^2 \, d\mu(z) + 2\mu(U)^{-1} \int_U |f(xz) - f(xy)|^2 \, d\mu(z).$$

En effectuant le changement de variable $z \to yz$ dans la deuxième intégrale, on obtient

$$|f(x) - f(xy)|^2 \leq 2\mu(U)^{-1} \int_U |f(x) - f(xz)|^2 \, d\mu(z) + 2\mu(U)^{-1} \int_{U^2} |f(xyz) - f(xy)|^2 \, d\mu(z).$$

Soit maintenant $x \in G$, $s \in G$. Il existe $s_0, \ldots, s_n \in U$ tel que $s = s_0 \cdots s_n$, où $n = |s|_U$ et $s_0 = e$, on a pour tout $i \in \{1, \ldots, n-1\}$, en utilisant l'inégalité ci-dessus appliquée à $xs_0 \cdots s_{i-1}$ au lieu de $x$ et $s_i$ au lieu de $y$ on obtient

$$|f(xs_0 \cdots s_{i-1}) - f(xs_0 \cdots s_i)|^2 \leq 2\mu(U)^{-1} \int_U |f(xs_0 \cdots s_{i-1}) - f(xs_0 \cdots s_{i-1}z)|^2 \, d\mu(z)$$
$$+ 2\mu(U)^{-1} \int_{U^2} |f(xs_0 \cdots s_{i-1}s_i) - f(xs_0 \cdots s_{i-1}s_i z)|^2 \, d\mu(z).$$

On en déduit

$$|f(xs_0 \cdots s_{i-1}) - f(xs_0 \cdots s_i)|^2 \leq \frac{2\mu(U)^{-1}}{r} \int_U |f(xs_0 \cdots s_{i-1}) - f(xs_0 \cdots s_{i-1}z)|^2 F_2(z) \, d\mu(z)$$
$$+ \frac{2\mu(U)^{-1}}{r} \int_{U^2} |f(xs_0 \cdots s_i) - f(xs_0 \cdots s_i z)|^2 F_2(z) \, d\mu(z),$$

donc en intégrant l'inégalité par rapport à $x$, on obtient

$$\int_G |f(xs_0 \cdots s_{i-1}) - f(xs_0 \cdots s_i)|^2 \, d\mu(x)$$
$$\leq \frac{2\mu(U)^{-1}}{r} \int_G \int_U |f(xs_0 \cdots s_{i-1}) - f(xs_0 \cdots s_{i-1}z)|^2 F_2(z) \, d\mu(z) \, d\mu(x)$$
$$+ \frac{2\mu(U)^{-1}}{r} \int_{U^2} \int_G |f(xs_0 \cdots s_i) - f(xs_0 \cdots s_i z)|^2 F_2(z) \, d\mu(z) \, d\mu(x)$$

et en effectuant des changements de variables dans les termes de droite, on obtient par la proposition 3

$$\int_G |f(xs_0 \cdots s_{i-1}) - f(xs_0 \cdots s_i)|^2 \, d\mu(x) \leq \frac{8\mu(U)^{-1}}{r} \varepsilon_{\nu_2}(f, f)$$

par ailleurs

$$|f(x) - f(xs)|^2 \leq n^2 \sum_{i=1}^n |f(xs_0 \cdots s_{i-1}) - f(xs_0 \cdots s_i)|^2$$

donc

$$\int_G |f(x) - f(xs)|^2 \, d\mu(x) \leq n^2 \frac{8\mu(U)^{-1}}{r} \varepsilon_{\nu_2}(f, f)$$

d'où

$$\int \int_G |f(x) - f(xs)|^2 F_1(s) \, d\mu(s) \, d\mu(x) \leq \int_G |s|_U^2 F_1(s) \, d\mu(s) \frac{8\mu(U)^{-1}}{r} \varepsilon_{\nu_2}(f, f). \qquad \Box$$

**Remarque.** *Dans les hypothèses de la proposition, l'existence d'un moment d'ordre 2 pour une densité $F$ ne dépend pas du choix de $U$, en effet; si $U$ et $U'$ sont deux ouverts relativement compacts voisinages de $e$ qui engendrent $G$ alors il existe un entier $p$ tel que $U' \subset U^p$ donc pour tout $x \in G$, $|x|_U \leq p|x|_{U'}$ par conséquent si $\int_G |x|_{U'}^2 F(x) \, dx$ est finie alors $\int_G |x|_U^2 F(x) \, dx$ est finie, d'où le résultat en échangeant les rôles de $U$ et $U'$.*

*De cette proposition et du résultat de comparaison des traces on déduit une proprieté d'invariance de la probabilité de retour.*

**Théorème 2.** *Soit $G$ un groupe localement compact séparable unimodulaire, compactement engendré, $F_1, F_2$ deux densités de probabilité symétriques bornées, tels qu'il existe des voisinages ouverts $U_1, U_2$ de $e$ relativement compacts, symétriques qui engendrent $G$ sur lesquels $\inf_{U_i} F_i > 0$ ; $i = 1, 2$, on suppose aussi que $F_i$ admet un moment d'ordre 2, c'est-à-dire que $\int_G |s|_{U_i}^2 F_i(s) \, d\mu(s)$ est finie, alors $F_1^{*(2t)}(e) \simeq F_2^{*(2t)}(e)$.*

**Démonstration.** Soit $\mathcal{A}$ la sous-algèbre hilbertienne formée des éléments bornés de $L^2(G)$, remarquons que

$$C_c(G) \subset \mathcal{A} \subset L^2(G),$$

où $C_c(G)$ est l'ensemble des fonctions continues sur $G$ et à valeurs complexes à support compact. Soit $R(\mathcal{A})$ l'adhérence faible de l'ensemble $\{R_f; f \in \mathcal{A}\}$. D'après ([5], 13.10.4, p. 272), $R(\mathcal{A})$ est une algèbre de von Neumann semi-finie. On définit une trace sur $R(\mathcal{A})^+$, en posant pour $S$ dans $R(\mathcal{A})^+$, $\tau(S) = (f|f)$ si $S^{1/2} = R_f$ avec $f$ borné et $\tau(S) = +\infty$ dans le cas contraire. Donc d'après ([6], théorème 1, p. 85) $\tau$ est une trace sur $R(\mathcal{A})$ fidèle, semi-finie et normale.

D'après la proposition 4, il existe une constante $C > 1$ telle que $\varepsilon_1 \leq C\varepsilon_2$, où $\varepsilon_i$ est la forme de Dirichlet associée à $F_i$, ce qui revient à dire que $I - R_{F_1} \leq C(I - R_{F_2})$, il s'ensuit par application du théorème 1 que $\tau(R_{F_2}^{2n}) \preceq \tau(R_{F_1}^{2n})$ par conséquent $F_2^{2n}(e) \preceq F_2^{2n}(e)$, d'où le résultat en échangeant les rôles de $F_1$ et $F_2$.

Le résultat du théorème 2 s'interprète grâce à la proposition 1 comme l'indépendance du comportement asymptotique de la probabilité de retour par rapport à la densité choisie. $\qquad \Box$

## 5. Marche aléatoire sur sol($K$)

Nous illustrons maintenant le théorème principal par l'exemple suivant. Notre reference principale sur les résultats et notations concernant les corps locaux est [16].

Soit $K$ un p-corps, on définit le groupe sol($K$) comme étant le produit semi direct de $K^*$ par $K^2$, où $K^*$ agit sur $K^2$ par les automorphismes $a \to \begin{pmatrix} a & 0 \\ 0 & a^{-1} \end{pmatrix}$, qui s'écrit alors

$$\text{sol}(K) = \left\{ \begin{pmatrix} a & 0 & x \\ 0 & a^{-1} & y \\ 0 & 0 & 1 \end{pmatrix}, a \in K^*, x, y \in K \right\}.$$

On désigne par $\text{mod}_K$ la fonction modulaire associée aux translations à gauche. On désigne respectivement par $R$ et $R^*$ la boule unité fermée et la boule unité ouverte associées à $\text{mod}_K$.

Soit $V = V_1^{-1} \cup V_1$ où

$$V_1 = \left\{ \begin{pmatrix} \pi^r u & 0 & R \\ 0 & \pi^{-r} u^{-1} & R \\ 0 & 0 & 1 \end{pmatrix}, r = 0, 1, -1, u \in R^* \right\}.$$

Ainsi $V$ est un compact symétrique, voisinage de $e$ qui engendre $\mathrm{sol}(K)$. On déduit alors que $\mathrm{Sol}(K)$ est localement compact, compactement engendré. Pour appliquer le théorème de stabilité on vérifie aisément que $\mathrm{Sol}(K)$ est de plus unimodulaire et résoluble donc moyennable.

Soit $\mu_1$ une mesure de Haar sur $K$. Comme $R$ est un compact de $K$, on peux choisir $\mu_1$ telle que $\mu_1(R) = 1$, de même $\mu_2$ mesure de Haar sur $K^*$ telle que $\mu_2(R^*) = 1$. Soit $\mu' = \mu_1 \otimes \mu_1 \otimes \mu_2$ la mesure produit sur $K^2 \times K^*$, et $\phi \colon K^2 \times K^* \to \mathrm{sol}(K)$, définie par $\phi(x, y, a) = \begin{pmatrix} a & 0 & x \\ 0 & a^{-1} & y \\ 0 & 0 & 1 \end{pmatrix}$. Considérons la mesure $\mu$ image de $\mu'$ par $\phi$, définie par $\mu(A) = \mu'(\phi^{-1}(A))$.

Il s'agit de montrer que $\mu$ est une mesure de Haar sur $\mathrm{sol}(K)$, on montre qu'elle est invariante par translation à gauche, pour cela il suffit d'utiliser un système de générateurs de la tribu borélienne de $\mathrm{sol}(K)$. Puisque les boréliens de la forme $S_1 \times S_2 \times S$ engendrent la tribu borélienne de $K \times K \times K^*$, par conséquent les boréliens de la forme

$$A = \phi(S_1 \times S_2 \times S)$$

où

$$S_1 \subset K, \qquad S_2 \subset K, \qquad S \subset K^*,$$

engendrent la tribu borélienne de $\mathrm{sol}(K)$. Soit $w = \phi(x_0, y_0, v)$, montrons alors que $\mu(wA) = \mu(A)$, où $A$ est choisi comme ci-dessus, en effet

$$\phi^{-1}(wA) = \{(vx + x_0, v^{-1}y + y_0, va) \mid (x, y) \in S_1 \times S_2, a \in S\} = S_1' \times S_2' \times S',$$

où

$$S_1' = vS_1 + x_0, \qquad S_2' = v^{-1}S_2 + y_0, \qquad S' = vS.$$

Par ailleurs $\mu_1(S_1') = \mu_1(vS_1) = \mathrm{mod}_K(v)\mu_1(S_1)$, $\mu_1(S_2') = \mu_1(v^{-1}S_2) = \mathrm{mod}_K(v^{-1})\mu_1(S_2)$ et $\mu_2(S') = \mu_2(vS) = \mu_2(S)$, donc

$$\mu'(S_1' \times S_2' \times S') = \mathrm{mod}_K(v)\mu_1(S_1)\mathrm{mod}_K(v^{-1})\mu_1(S_2)\mu_2(S) = \mu'(S_1 \times S_2 \times S)$$

ce qui entraîne $\mu(wA) = \mu'(S_1' \times S_2' \times S') = \mu'(S_1 \times S_2 \times S) = \mu(A)$.

En conclusion nous obtenons une mesure de Haar à gauche sur $\mathrm{sol}(K)$.

On considère pour tout entier $n > 0$,

$$\Omega_n = \left\{ \begin{pmatrix} a & 0 & x \\ 0 & a^{-1} & y \\ 0 & 0 & 1 \end{pmatrix}, a \in K^*, \mathrm{mod}_K(a) \in [q^{-n}, q^n], x, y \in \pi^{-n}R \right\}.$$

En appliquant des arguments analogues à ceux de [4] on obtiendra le comportement de la marche sur $\mathrm{sol}(K)$, nous commençons par calculer dans cette partie le volume de ces boîtes. On a $\mu(\Omega_n) = \mu_1(\pi^{-n}R)^2 \mu_2(B)$, où

$$B = \{a \in K^*, \mathrm{mod}_K(a) \in [q^{-n}, q^n]\}$$

et on a

$$\mu_1(\pi^{-n}R) = \mathrm{mod}_K(\pi^{-n})\mu_1(R) = q^n$$

et d'après le théorème 3.2, $\Gamma = \mathrm{mod}_K(K^*) = \{q^m, m \in \mathbb{Z}\}$ donc

$$B = \bigcup_{k=-n}^{n} \{a \in K^* \mid \mathrm{mod}_K(a) = q^k\} = \bigcup_{k=-n}^{n} \{a \in K^* \mid \mathrm{mod}_K(\pi^k a) = 1\} = \bigcup_{k=-n}^{n} \pi^{-k} R^*$$

donc

$$\mu_2(B) = (2n+1)\mu_2(R^*) = 2n+1$$

d'où $\mu(\Omega_n) = (2n+1)q^{2n}$.

On définit la marche sur $\mathrm{sol}(K)$ par la donnée de la densité $F = \frac{1}{\mu(V)} 1_V$.

Pour appliquer le théorème 1, il nous reste à prouver que le groupe est moyennable, pour cela on considère les flèches suivantes

$$d: \mathrm{sol}(K) \to K^*, \qquad \begin{pmatrix} a & 0 & x \\ 0 & a^{-1} & y \\ 0 & 0 & 1 \end{pmatrix} \mapsto a$$

et

$$w: K^* \to \mathbb{Z}, \qquad a \mapsto -\log_q(\mathrm{mod}_K(a))$$

on a alors la suite exacte

$$0 \to K \times K \to \mathrm{sol}(K) \to \mathbb{Z} \to 0$$

par conséquent d'après le théorème 1.2, $\mathrm{sol}(K)$ est moyennable.

**Lemme 3.** *Si $s_1, \ldots, s_t$ sont des éléments de $V$ tels que pour tout $i \in \{1, \ldots, t\}$, $\mathrm{mod}_K(d(s_1, \ldots, s_i)) \in [q^{-n}, q^n]$ alors $s_1, \ldots, s_t \in \Omega_n$.*

**Démonstration.** On raisonne par récurrence sur $t$, pour $t=1$, on a $s_1 \in V$ donc $d(s_1)$ s'écrit $d(s_1) = \pi^r u$, où $u \in R^*, r \in \{0, 1, -1\}$ donc $\mathrm{mod}_K(d(s_1)) = q^{-r} \in [q^{-1}, q]$ donc $s_1 \in \Omega_n$.

Supposons la propriété pour $t$ et considérons $s_1, \ldots, s_{t+1} \in V$ tels que l'on a pour tout $i \in \{1, \ldots, t+1\}$; $\mathrm{mod}_K(d(s_1 \cdots s_i)) \in [q^{-n}, q^n]$, en particulier $s_1 \cdots s_t \in \Omega_n$ donc s'écrit

$$s_1 \cdots s_t = \begin{pmatrix} a & 0 & x \\ 0 & a^{-1} & y \\ 0 & 0 & 1 \end{pmatrix}, \quad \mathrm{mod}_K(a) \in [q^{-n}, q^n], x, y \in \pi^{-n} R.$$

Si $s_{t+1} \in V_1$, alors on peut écrire $s_{t+1} = \begin{pmatrix} \pi^r u & 0 & x_0 \\ 0 & \pi^{-r} u^{-1} & y_0 \\ 0 & 0 & 1 \end{pmatrix}$, où $u \in R^\times$, $r \in \{0, 1, -1\}$, $x_0, y_0 \in R$ donc

$$s_1 \cdots s_t s_{t+1} = \begin{pmatrix} \pi^r u a & 0 & ax_0 + x \\ 0 & \pi^{-r} u^{-1} a^{-1} & a^{-1} y_0 + y \\ 0 & 0 & 1 \end{pmatrix}.$$

Par ailleurs $\mathrm{mod}_K(d(s_1 \cdots s_{t+1})) \in [q^{-n}, q^n]$ par hypothèse du lemme, et on a

$$\mathrm{mod}_K(ax_0 + x) \leq \max(\mathrm{mod}_K(ax_0), \mathrm{mod}_K(x)) \leq q^n$$

donc $ax_0 + x \in \pi^{-n} R$, de même $a^{-1} y_0 + y \in \pi^{-n} R$, d'où $s_1 \cdots s_{t+1} \in \Omega_n$.

Si $s_{t+1} \in V_1^{-1}$, alors on peut écrire $s_{t+1} = \begin{pmatrix} \pi^{-r}u^{-1} & 0 & \pi^{-r}u^{-1}x_0 \\ 0 & \pi^r u & \pi^r u y_0 \\ 0 & 0 & 1 \end{pmatrix}$, où $u \in R^\times$, $r \in \{0, 1, -1\}$, $x_0, y_0 \in R$, donc

$$s_1 \cdots s_t s_{t+1} = \begin{pmatrix} \pi^{-r}u^{-1}a & 0 & \pi^{-r}u^{-1}ax_0 + x \\ 0 & \pi^r u a^{-1} & \pi^r u a^{-1} y_0 + y \\ 0 & 0 & 1 \end{pmatrix}.$$

On a

$$\mathrm{mod}_K(\pi^{-r}u^{-1}ax_0 + x) \leq \max(\mathrm{mod}_K(\pi^{-r}u^{-1}ax_0), \mathrm{mod}_K(x)) \leq \max(\mathrm{mod}_K(\pi^{-r}a), q^n).$$

Par ailleurs $\mathrm{mod}_K(d(s_1 \cdots s_{t+1})) \leq q^n$ donc $\mathrm{mod}_K(\pi^{-r}u^{-1}ax_0 + x) \leq q^n$, par conséquent $\pi^{-r}u^{-1}ax_0 + x \in \pi^{-n}R$, et on montre de même que $\pi^r u a^{-1} y_0 + y \in \pi^{-n}R$, par suite $s_1 \cdots s_t s_{t+1} \in \Omega_n$. □

En utilisant le lemme 3, on obtient le théorème.

**Théorème 3.** *La probabilité de retour de la marche aléatoire sur* $\mathrm{sol}(K)$ *à le comportement asymptotique de* $\exp(-t^{1/3})$.

**Démonstration (Comparer avec [1]).** Soit $S_1, S_2, \ldots, S_t$ les variables aléatoires obtenues en projetant $Z_1, \ldots, Z_t$ sur $\mathbb{Z}$, c'est a dire $S_i = \omega \circ d \circ Z_i$. $(S_i)$ est la marche aléatoire sur $\mathbb{Z}$ définie par la mesure $\frac{1}{5}(2\delta_{-1} + \delta_0 + 2\delta_1)$.

En utilisant le lemme 3 on obtient

$$P(Z_{2t} \in \Omega_n) \geq P(Z_1 \in \Omega_n, Z_2 \in \Omega_n, \ldots, Z_{2t} \in \Omega_n) \geq P(S_1 \in [-n, n], \ldots, S_{2t} \in [-n, n])$$

et d'après [10], on a

$$P(S_1 \in [-n, n], \ldots, S_{2t} \in [-n, n]) \succeq \exp\left(-\frac{t}{n^2}\right)$$

d'où

$$P(Z_{2t} \in \Omega_n) \geq P(Z_1 \in \Omega_n, Z_2 \in \Omega_n, \ldots, Z_{2t} \in \Omega_n) \succeq \exp\left(-\frac{t}{n^2}\right).$$

Par ailleurs, en utilisant l'inégalité de Cauchy–Schwartz

$$P(Z_{2t} \in \Omega_n) = \int_{\Omega_n} F^{*(2t)}(g) \, d\mu(g) \leq \mu(\Omega_n) F^{*(2t)}(e)$$

par conséquent

$$F^{*(2t)}(e) \succeq \exp\left(-\frac{t}{n^2}\right)/q^{2n} \succeq \exp\left(-\frac{t}{n^2} - n\right),$$

et en choisissant $n$ tel que $n^3 \simeq t$ on obtient $F^{*(2t)}(e) \succeq \exp(-t^{1/3})$.

On a d'après ([8], théorème 4.1, p. 683), pour tout groupe de croissance exponentielle, $F^{*(2t)}(e) \preceq \exp(-t^{1/3})$. On en déduit que sur $\mathrm{sol}(K)$, $F^{*(2t)}(e) \simeq \exp(-t^{1/3})$.

En conclusion pour toute densité $F$ sur $\mathrm{sol}(K)$ vérifiant les hypothèses du théorème 2.3, on a

$$P(Z_{2n} \in \Omega) \simeq \exp(-n^{1/3}). \qquad \square$$

## 6. Conclusions et perspectives

Soit $G$ un groupe localement compact compactement engendré unimodulaire. Le théorème principal montre qu'il existe un comportement asymptotique naturellement associé à $G$; la décroissance de $F^{*(2n)}(e)$.

Quels sont les incidences des propriétés algébriques de $G$ sur la décroissance de $F^{*(2n)}(e)$?

Nous pensons que le cas de $\text{sol}(K)$ traité en détail est représentatif pour la classe des groupes algébriques résolubles unimodulaire sur un corps local.

Une perspective de recherche que nous n'avons malheureusement pas eu le temps de concrétiser et formaliser dans ce texte est celle de la comparaison entre probabilité de retour dans un groupe et probabilité de retour dans un de ses sous groupes. Cette comparaison permet par exemple d'expliquer de manière assez conceptuelle pourquoi la probabilité de retour dans le groupe $\mathbb{Z}/2\mathbb{Z} \wr \mathbb{Z}$ se comporte en $\exp(-t^{1/3})$. Ce groupe est un réseau cocompact dans le groupe $\text{sol}(K)$ où $K$ est le corps local $\mathbb{F}_2((t))$ des séries formelles à coefficients dans le corps à deux éléments.